\documentclass[a4paper,11pt]{article}
\usepackage{enumerate}
\usepackage{amsmath,amsthm,amssymb}

\title{Traces Without Maximal Chains}
\author{Ta Sheng Tan \thanks{Department of Pure Mathematics and Mathematical Statistics, Centre for Mathematical Sciences, University of Cambridge, Wilberforce Road, Cambridge CB3 0WB, United Kingdom. Email: T.S.Tan@dpmms.cam.ac.uk.}}
\newtheorem{thm}{Theorem}[section]
\newtheorem{lemma}[thm]{Lemma}
\newtheorem{conj}[thm]{Conjecture}
\theoremstyle{remark} \newtheorem{claim}{Claim}
\theoremstyle{remark} \newtheorem{pfclaim}{Proof of Claim}

\begin{document}

\maketitle

\begin{abstract}
The trace of a family of sets $\mathcal{A}$ on a set $X$ is $\mathcal{A}|_X=\{A\cap X:A\in \mathcal{A}\}$. 
If $\mathcal{A}$ is a family of $k$-sets from an $n$-set such that for any $r$-subset $X$ the trace $\mathcal{A}|_X$ does not contain a maximal chain, then how large can $\mathcal{A}$ be? 
Patk\'os conjectured that, for $n$ sufficiently large, the size of $\mathcal{A}$ is at most $\binom{n-k+r-1}{r-1}$. Our aim in this paper is to prove this conjecture.
\end{abstract}

\section{Introduction}
Let $[n]$ denote the set of integers $\{1,2,\ldots,n \}$. Given a set $X$ we write $\mathcal{P}(X)$ for its power set and $X^{(k)}$ for the set of all its $k$-element subsets (or $k$-subsets). The \textit{trace} of a family $\mathcal{A}$ of sets on a set $X$ is $\mathcal{A}|_X=\{A\cap X:A\in \mathcal{A}\}$.
\\\\
Vapnik and Chervonenkis \cite{vapnik}, Sauer \cite{sauer} and Shelah \cite{shelah}  independently showed that if $\mathcal{A} \subset \mathcal{P}([n])$ is a family with more than $\sum_{i=0}^{k-1} \binom{n}{i}$ sets, then there is a $k$-subset $X$ of $[n]$ such that $\mathcal{A}|_X=\mathcal{P}(X)$. This bound is sharp, as shown for example by the family $\{A\in [n]:|A|< k\}$, but no characterisation for the extremal families is known.
\\\\
The uniform case of the problem was considered by Frankl and Pach \cite{frankl}. They proved that if $\mathcal{A} \subset [n]^{(k)}$ is a family with more than $\binom{n}{k-1}$ sets, then there is a $k$-subset $X$ of $[n]$ such that $\mathcal{A}|_X=\mathcal{P}(X)$. This bound is not sharp and was improved later by Mubayi and Zhao \cite{mubayi}, but the exact bound is still unknown.
\\\\
While the above problems concern families with traces not containing the power set, Patk\'os \cite{patkos1,patkos2} considered the case of families with traces not containing a maximal chain. Here a \textit{maximal chain} of a set $X$ is a family of the form $X_0 \subset X_1 \subset X_2 \ldots \subset X_r=X$, where $|X_i|=i$ for all $i$. He proved in \cite{patkos2} that if $\mathcal{A} \subset \mathcal{P}([n])$ is a family with more than $\sum_{i=0}^{k-1} \binom{n}{i}$ sets, then there is a $k$-subset $X$ of $[n]$ such that the trace $\mathcal{A}|_X$ contains a maximal chain of $X$, with the only extremal families being $\{A\in [n]:|A|< k\}$ and $\{A\in [n]:|A|> n-k\}$. This beautiful result is an extension of the result of Vapnik and Chervonenkis, Sauer and Shelah. For the $k$-uniform case, he proved in \cite{patkos1} that $\{A \in [n]^{(k)}:1\in A\}$ is an extremal family for $n$ sufficiently large: in other words, if $\mathcal{A}\subset [n]^{(k)}$ has more than $\binom{n-1}{k-1}$ sets, then there is a $k$-subset $X$ of $[n]$ such that the trace $\mathcal{A}|_X$ contains a maximal chain of $X$. He also proved the stability of this extremal family. He further conjectured that for any $k\geq r\geq 2$, if $n$ is sufficiently large and $\mathcal{A} \subset [n]^{(k)}$ has more than $\binom{n-k+r-1}{r-1}$ sets there is an $r$-subset $X$ of $[n]$ such that the trace $\mathcal{A}|_X$ contains a maximal chain of $X$.
\\\\
In this paper, we prove this conjecture. Our proof also shows that the only extremal families are of the form $\{A \in [n]^{(k)}:D\subset A\}$, for some $(k-r+1)$-subset $D$ of $[n]$. 
\\\\
For $n\geq k$ and $r\geq 1$, we define $W(n,k,r)$ to be the maximum size of a $k$-uniform family $\mathcal{A}\subset [n]^{(k)}$ with the property that for any $r$-subset $X$ the trace $\mathcal{A}|_X$ does not contain a maximal chain of $X$. Thus, our main result is to show that for $k\geq r$, $W(n,k,r)=\binom{n-k+r-1}{r-1}$, provided $n$ is sufficiently large.
\\\\
Patk\'os \cite{patkos1} proved the case $k=r$ using a stability theorem of Hilton and Milner \cite{hilton} about intersecting families. Our proof for the general case, which does not use Patk\'os' result, is self-contained, and in fact also yields a simpler proof of Patk\'os' result.

\section{Main Result}

The idea of the proof is as follows. We split the problem into two cases: the case when $\mathcal{A}$ is intersecting and the case when $\mathcal{A}$ is non-intersecting. It turns out that the former case can be done in a straightforward way by induction. For the latter case, we reduce the problem to considering extremal families with traces not containing an ``almost'' maximal chain. Here, an \textit{almost maximal chain} of a set $X$ is a maximal chain of $X$ without the empty set, i.e. a family of the form $\{X_1 \subset X_2 \ldots \subset X_r=X: |X_i|=i\}$. (Interestingly, almost maximal chains were also considered by Patk\'os \cite{patkos1}.) At this point, it looks like we might need to further reduce the problem to considering extremal families with traces not containing an ``almost almost'' maximal chain and so on. But luckily, this is not the case, as one can bound the sizes of extremal families with traces not containing an almost maximal chain in terms of the sizes of extremal families with traces not containing a maximal chain.
\\\\
For $n\geq k$ and $r\geq 2$, we define $U(n,k,r)$ for the maximum size of a $k$-uniform family $\mathcal{A}\subset [n]^{(k)}$ with the property that for any $r$-subset $X$ the trace $\mathcal{A}|_X$ does not contain an almost maximal chain of $X$. With this notation, we have the following lemma.

\begin{lemma}
For $k,r\geq2$, $U(n,k,r)\leq \frac{n}{k} W(n-1,k-1,r-1)$. 
\end{lemma}

\begin{proof}
 Let $\mathcal{A}\subset [n]^{(k)}$ be such that for any $r$-subset $X$ the trace $\mathcal{A}|_{X}$ does not contain an almost maximal chain of $X$.
\\\\
For each $x \in [n]$, define $\mathcal{B}_{\{x\}}=\{A \in \mathcal{A}: x \in A\}$. We then claim that $|\mathcal{B}_{\{x\}}|\leq W(n-1,k-1,r-1)$.
Suppose not, then the family $\mathcal{C}_{\{x\}}=\{B\setminus \{x\}:B \in \mathcal{B}_{\{x\}}\}$ is a $(k-1)$-uniform family in $([n]\setminus \{x\})^{(k-1)}$ with size greater than $W(n-1,k-1,r-1)$. By definition, there exists an $(r-1)$-subset $X$ not containing $x$ such that the trace $(\mathcal{C}_{\{x\}})|_X$ contains a maximal chain of $X$. So, $(\mathcal{B}_{\{x\}})|_{X\cup \{x\}}$ contains an almost maximal chain of $X\cup \{x\}$. This is a contradiction.
\\\\
By averaging over all possible $x$, we have $|\mathcal{A}| \leq \frac{n}{k} W(n-1,k-1,r-1).$
\end{proof}

\noindent
We can now prove our main theorem. Note that for $k<r$, we have $W(n,k,r)=\binom{n}{k}$. Also, $W(n,k,1)=1$.

\begin{thm} \label{maintheorem}
Let $k \geq r-1$. Then there exists an $n_0(k,r)$ such that for any $n \geq n_0(k,r)$, $W(n,k,r) = \binom{n-k+r-1}{r-1}$.
\end{thm}

\begin{proof}
We use induction on $r$, and for fixed $r$ induction on $k$. The theorem is clearly true for $r=1$. So fix $r > 1$ and suppose that the theorem is true for $r-1$ (and all $k\geq r-2$). For our given value of $r$, the theorem is trivially true for $k=r-1$.
\\\\
Now fix $k \geq r$ and suppose that the theorem is true for $k-1$. Let $\mathcal{A} \subset [n]^{(k)}$ be a $k$-uniform family such that for any $r$-subset $X$ the trace $\mathcal{A}|_X$ does not contain a maximal chain of $X$.
\\\\
\textbf{Case 1: $\mathcal{A}$ is intersecting.}

\noindent
We may assume $\bigcap_{A\in \mathcal{A}} A =\emptyset$. Otherwise, let $x \in \bigcap_{A\in \mathcal{A}} A$, and then by induction, we have 
\begin{align*}
  |\mathcal{A}|=|\{A\setminus {x}: A \in \mathcal{A}\}|&\leq W(n-1,k-1,r)\\
						       &=\binom{n-k+r-1}{r-1}, \mbox{ as required.}
\end{align*}

\noindent
Now let $ l=\min\{|A\cap B|:A,B \in \mathcal{A}\}$: so $l\geq 1$.
Pick $A,B$ such that $|A\cap B|=l.$ We may then write $\mathcal{A}=\mathcal{A}_1 \cup \mathcal{A}_2$, where
$\mathcal{A}_1=\{C\in \mathcal{A}:C\supset A\cap B\}$ and $\mathcal{A}_2=\mathcal{A} \setminus \mathcal{A}_1$.
Since $\bigcap_{A\in \mathcal{A}} A = \emptyset$, we have $\mathcal{A}_2 \neq \emptyset$. Pick $D\in \mathcal{A}_2$. Note that $(A\cap B)\setminus D \neq \emptyset$.

\begin{claim}
 $|\mathcal{A}_1|\leq 8^k \binom{n-k}{r-2}$.
\end{claim}
\begin{pfclaim}
 For each $S \subset A\cup B\cup D$, define $\mathcal{B}_S = \{C\in \mathcal{A}_1:C\cap (A \cup B\cup D)=S\}$ and $\mathcal{C}_S = \{F\setminus S: F\in \mathcal{B}_S\}$. $\mathcal{B}_S$ is non-empty only if $S\supset A\cap B$. Suppose $|\mathcal{B}_S|>W(n-|A\cup B\cup D|,k-|S|,r-1)$, then there exists an $(r-1)$-subset $X \subset [n]\setminus (A\cup B\cup D)$ such that the trace $\mathcal{C}_S|_X$ contains a maximal chain of $X$. Pick $a\in (A\cap B)\setminus D$, then $\mathcal{B}_S|_{X\cup \{a\}}$ contains an almost maximal chain of $X\cup \{a\}$ and $D\cap (X\cup \{a\})=\emptyset$. This is a contradiction as $\mathcal{A}|_{X\cup \{a\}}$ would contain a maximal chain of $X\cup \{a\}$. Hence, $|\mathcal{B}_S|\leq W(n-k,k-|S|,r-1)\leq \binom{n-k}{r-2}$. This completes the proof of the claim.
\end{pfclaim}

\begin{claim}
 $|\mathcal{A}_2|\leq 4^k \binom{n-k}{r-2}$.
\end{claim}
\begin{pfclaim}
 As before, for each $S\subset (A\cup B)$, we define $\mathcal{B}_S = \{C\in \mathcal{A}_2:C\cap (A \cup B)=S\}$ and $\mathcal{C}_S = \{F\setminus S: F\in \mathcal{B}_S\}$. By the minimality of $l$, $\mathcal{B}_S$ is non-empty only if $S\cap (A\setminus B) \neq \emptyset$. Suppose $|\mathcal{B}_S|>W(n-|A\cup B|,k-|S|,r-1)$, then there exists an $(r-1)$-subset $X \subset [n]\setminus (A\cup B)$ such that the trace $\mathcal{C}_S|_X$ contains a maximal chain of $X$. Pick $a\in S\cap (A\setminus B)$, then $\mathcal{B}_S|_{X\cup \{a\}}$ contains an almost maximal chain of $X\cup \{a\}$ and $B\cap (X\cup \{a\})=\emptyset$. This is a contradiction as $\mathcal{A}|_{X\cup \{a\}}$ would contain a maximal chain of $X\cup \{a\}$. Hence, $|\mathcal{B}_S|\leq W(n-k,k-|S|,r-1)\leq \binom{n-k}{r-2}$. This completes the proof of the claim.
\end{pfclaim}

\noindent
So $|\mathcal{A}|=|\mathcal{A}_1|+|\mathcal{A}_2|\leq (8^k+4^k)\binom{n-k}{r-2}$, which is certainly at most $\binom{n-k+r-1}{r-1}$ for $n$ sufficiently large.
\\\\
\textbf{Case 2: $\mathcal{A}$ is non-intersecting.}

\noindent
Let $A$ and $B$ be in $\mathcal{A}$ such that $A\cap B=\emptyset$. We may then write $\mathcal{A}=\mathcal{A}_1 \cup \mathcal{A}_2$, where $\mathcal{A}_1=\{C\in \mathcal{A}:C\cap A \neq \emptyset\}$ and $\mathcal{A}_2=\mathcal{A} \setminus \mathcal{A}_1$. It is easy to see that $\mathcal{A}_2$ is a $k$-uniform family in $([n]\setminus A)^{(k)}$ such that for any $r$-subset $X$ in $[n]\setminus A$ the trace $\mathcal{A}_2|_X $ does not contain an almost maximal chain of $X$. So, 
\begin{align*}
  |\mathcal{A}_2| &\leq U(n-k,k,r) \\
		  &\leq \frac{n-k}{k} W(n-k-1,k-1,r-1)\\
	          &\leq \frac{n-k}{k} W(n,k-1,r-1)\\
	          &= \frac{n-k}{k}\binom{n-k+r-1}{r-2}\\
	          &= \frac{r-1}{k}\binom{n-k+r-1}{r-1}.
\end{align*}

\noindent
We are now left to bound the size of $\mathcal{A}_1$. 

\begin{claim}
 $|\mathcal{A}_1|\leq 4^k \binom{n-k}{r-2}$.
\end{claim}
\begin{pfclaim}
 Again, for each $S\subset (A\cup B)$, we define $\mathcal{B}_S = \{C\in \mathcal{A}_1:C\cap (A \cup B)=S\}$ and $\mathcal{C}_S = \{F\setminus S: F\in \mathcal{B}_S\}$. $\mathcal{B}_S$ is non-empty only if $S\cap A \neq \emptyset$. Suppose $|\mathcal{B}_S|>W(n-|A\cup B|,k-|S|,r-1)$, then there exists an $(r-1)$-subset $X \subset [n]\setminus (A\cup B)$ such that the trace $\mathcal{C}_S|_X$ contains a maximal chain of $X$. Pick $a\in S\cap A$, then $\mathcal{B}_S|_{X\cup \{a\}}$ contains an almost maximal chain of $X\cup \{a\}$ and $B\cap (X\cup \{a\})=\emptyset$. This is a contradiction as $\mathcal{A}|_{X\cup \{a\}}$ would contain a maximal chain of $ X\cup \{a\}$. Hence, $|\mathcal{B}_S|\leq W(n-k,k-|S|,r-1)\leq \binom{n-k}{r-2}$. This completes the proof of the claim.
\end{pfclaim}

\noindent
So we have
\begin{align*}
 |\mathcal{A}|&=|\mathcal{A}_1|+|\mathcal{A}_2|\\
	      &\leq 4^k \binom{n-k}{r-2} + \frac{r-1}{k}\binom{n-k+r-1}{r-1}.
\end{align*}

\noindent
As $k>r-1$, this is certainly at most $\binom{n-k+r-1}{r-1}$ for $n$ sufficiently large. 
\end{proof}

\noindent
Note that for a fixed $r$, equality can only hold (for $n$ sufficiently large) if $\bigcap_{A\in \mathcal{A}}A\neq \emptyset$ for each of the induction steps. This shows that the only extremal families are of the form $\{A \in [n]^{(k)}:D\subset A\}$, for some $(k-r+1)$-subset $D$ of $[n]$. 

\section{Remarks}

In this section, we give a few remarks relating to the proof of Theorem~\ref{maintheorem}.
\\\\
To give an explicit $n_0(k,r)$, we need $\binom{n-k+r-1}{r-1} \geq \max\{4^k \binom{n-k}{r-2} + \frac{r-1}{k}\binom{n-k+r-1}{r-1},(8^k+4^k)\binom{n-k}{r-2}\}$ and so we can take $n_0(k,r)=r8^k$. This is clearly not optimal. A more careful case analysis shows that $n_0(k,2) = 2k$ and trivially $n_0(k,1)=k$. This suggests that $n_0(k,r)=rk$ might suffice, but actually we believe that  $n_0(k,r)$ can be as small as $2k+1$.

\begin{conj}
 For $k\geq r\geq 3$, Theorem~\ref{maintheorem} holds with $n_0(k,r)=2k+1$.
\end{conj}

\noindent
Note that this cannot be improved to $2k$ in general - for example, one can check that $n_0(3,3)=7$.
\\\\
While we have shown that there is a unique (up to permutation of the ground set) extremal family for $n$ large, we are also interested in finding extremal families for \textit{all} $n\geq k$. For the case $r=2$ and $k+1 \leq n < 2k$, $[k+1]^{(k)}$ is the unique (up to permutation) extremal family.

\begin{conj}
Let $r\geq 2$. For $k+r-1 \leq n < 2k$, $W(n,k,r)=\binom{k+r-1}{k}$ and the only extremal family is of the form $[k+r-1]^{(k)}$. For $n>2k$, $W(n,k,r)=\binom{n-k+r-1}{r-1}$ and the only extremal family is of the form $\{A \in [n]^{(k)}:1,2,\ldots,k-r+1\in A\}$.
\end{conj}

\subsection*{Acknowledgement}
 
The author is very thankful to Allan Siu Lun Lo for helpful discussions, and to Imre Leader for his invaluable comments.

\end{document}